\documentclass[conference, final]{IEEEtran}
\IEEEoverridecommandlockouts
\usepackage{graphicx}
\usepackage{textcomp}
\usepackage{xcolor}

\usepackage{algorithm}
\usepackage{algpseudocode}
\usepackage{bbm}
\usepackage{enumerate}
\usepackage{multirow}
\usepackage{amsmath,amssymb,amsfonts}
\usepackage{subfig}
\usepackage{fancyhdr}

\newcommand{\real}{\ensuremath{\mathbb{R}}}

\newcommand{\abs}[1]{\left|#1\right|}
\newcommand{\norm}[1]{\left\lVert#1\right\rVert}

\newcommand{\E}{\ensuremath{\mathbb{E}}}

\def\BibTeX{{\rm B\kern-.05em{\sc i\kern-.025em b}\kern-.08em
    T\kern-.1667em\lower.7ex\hbox{E}\kern-.125emX}}
    
\renewcommand{\baselinestretch}{1.0}

\begin{document}

\title{Chance Constraint Tuning for Optimal Power Flow }

\author{\IEEEauthorblockN{Ashley M. Hou and Line A. Roald}
\IEEEauthorblockA{\textit{Electrical and Computer Engineering,
University of Wisconsin-Madison, Madison, USA} \\
\{amhou, roald\}@wisc.edu}
\vspace{-0.2in}
\thanks{This research is supported by the Department of Energy, Office of Science, Office of Advanced Scientific Computing Research, Applied Mathematics program under Contract Number DE-AC02-06CH11347.}
}

\maketitle

\begin{abstract}
In this paper, we consider a chance-constrained formulation of the optimal power flow problem to handle uncertainties resulting from renewable generation and load variability. We propose a tuning method that iterates between solving an approximated reformulation of the optimization problem and using a posteriori sample-based evaluations to refine the reformulation. Our method is applicable to both single and joint chance constraints and does not rely on any distributional assumptions on the uncertainty. In a case study for the IEEE 24-bus system, we demonstrate that our method is computationally efficient and enforces chance constraints without over-conservatism.
\end{abstract}

\section{Introduction}
\bstctlcite{}
The optimal power flow (OPF) problem aims to minimize total generation costs while enforcing physical system constraints on power flow balance and generator and line limits. However, increased penetration of renewable energy production has introduced uncertainties than can render solutions obtained with traditional deterministic methods insecure. To handle these variabilities, several stochastic versions of the OPF have been proposed. We consider a chance-constrained OPF (CC-OPF) formulation, where constraints are required to be satisfied with probability greater than $1 - \epsilon$, with $\epsilon$ denoting the acceptable violation probability. 
Chance constraints are an intuitive way to account for uncertainty, and are used in practice for, e.g., reserve dimensioning \cite{abbaspourtorbati2016}. However, chance-constrained problems are generally difficult to solve to optimality. Some special cases (e.g., Gaussian uncertainty) admit exact analytical reformulations \cite{roaldAnalyticalReformulationSecurity2013, bienstock2014}, leading to tractable problem formulations at the expense of strong limitations on the type of distribution. Other methods such as the scenario approach \cite{calafioreScenarioApproachRobust2006, vrakopoulou2013} or distributionally robust methods \cite{summers2014, roaldSecurityConstrainedOptimal2015, lubin2015, li2019distributionally} require less stringent assumptions on the distribution, and achieve feasible solutions by prioritizing chance constraint satisfaction and problem tractability over optimality. In many cases, this leads to very conservative (i.e., sub-optimal) solutions or infeasibility. Sample average approximations, which approximate the probabilistic constraint based on a large set of samples, can find an optimal solution by identifying a subset of samples where violations are allowed. These problems are typically formulated as mixed-integer problems \cite{luedtkeSampleApproximationApproach}, making numerical tractability a challenge. However, recent work has investigated continuous approximations, e.g. \cite{pena2019dc}. 

\begin{figure}
	\centering
	\includegraphics[width=0.95\linewidth]{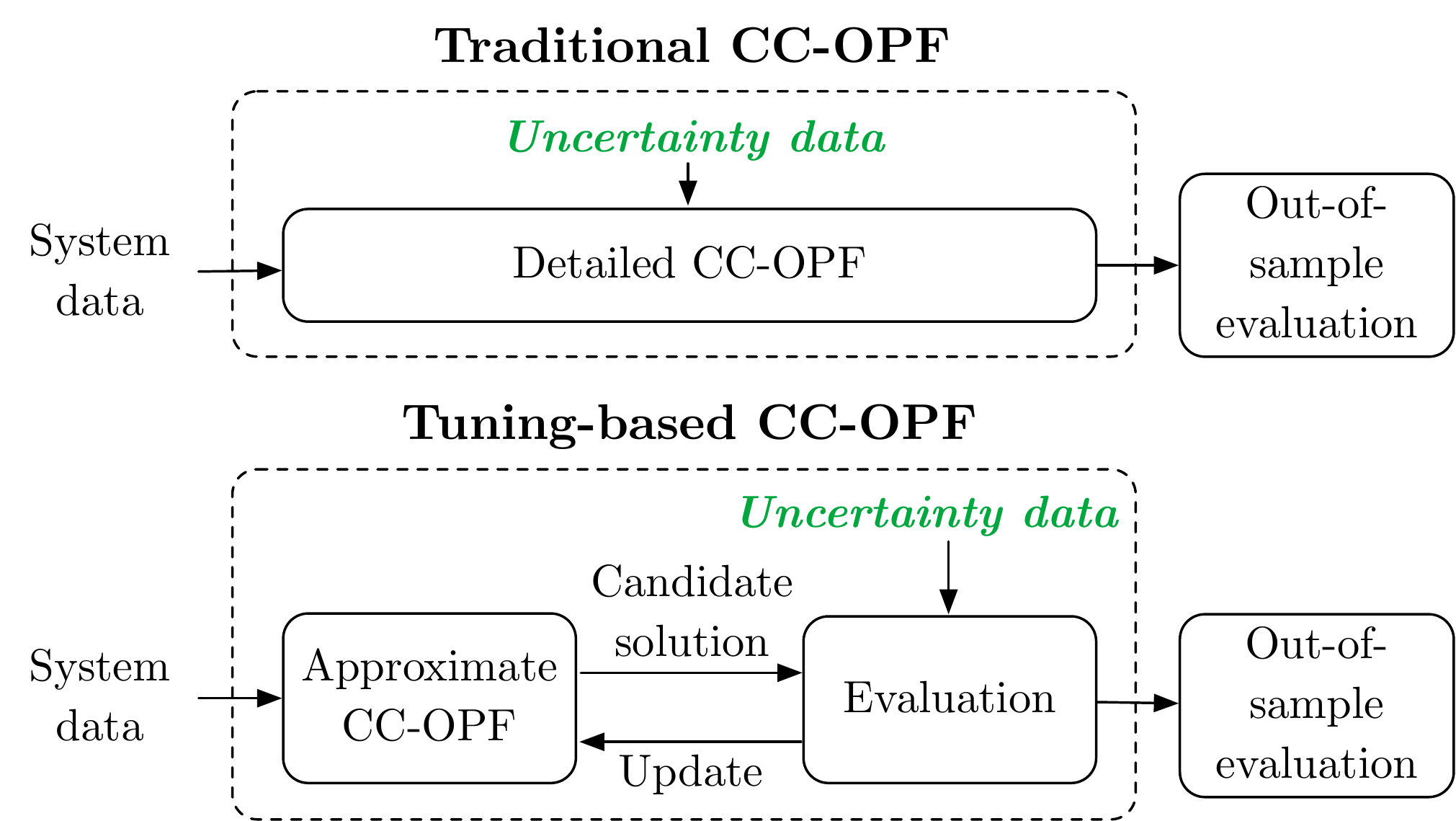}
	\caption{\small{Traditional solution methods for CC-OPF (top) include uncertainty information directly in a detailed problem formulation. The tuning method (bottom) uses a simple, approximate problem formulation which is iteratively updated using uncertainty data.}}
	\vspace{-0.25in}
	\label{fig:trad-block-diagram}
\end{figure}

As illustrated in Figure \ref{fig:trad-block-diagram} (top), these existing methods all include information about the uncertainty distribution directly in the optimization problem formulation, either in the form of parameters (i.e., the mean and covariance matrix) or through samples. Typically, there is a trade-off between the computational tractability (which is reduced as more information is taken into account) and the quality of the resulting solution (which is improved with more information).

In this paper, we propose utilizing results from a posteriori sample-based tests to improve a simple CC-OPF formulation. Specifically, we use data-driven parameter tuning to achieve a desired violation probability for both joint and single chance constraints. As illustrated in Figure \ref{fig:trad-block-diagram} (bottom), 
our method iterates between two steps: (i) solving an approximate chance-constrained problem to obtain a candidate solution and (ii) using samples to evaluate the candidate solution and update the problem approximation. In this way, we can utilize information from samples without including them in the optimization problem itself. The approximate optimization problem is computationally inexpensive and easily solved with commercial solvers. Therefore, although the optimization problem must be solved at each iteration of the tuning process, our method may be less computationally intense, yet more accurate, than alternative methods. Moreover, our method is applicable to both single and joint chance-constrained problems.

We note that our method is not the first to consider iterative tuning \cite{lam2019combating, mezghani2019stochastic, pena2019dc} and online updates \cite{oldewurtelAdaptivelyConstrainedStochastic2013, oldewurtelAdaptivelyConstrainedStochastic2015} to improve solutions to chance-constrained problems. For example, \cite{lam2019combating} addresses the theory of tuning safety parameters for a generic chance-constrained problem while still maintaining probabilistic guarantees in a limited data regime.
In \cite{pena2019dc}, an iterative tuning process was proposed to more accurately determine the value of a safety parameter in joint chance-constrained OPF. The tuning process in \cite{pena2019dc} can be characterized as \emph{fine-tuning}, where smaller adjustments are made to an accurate model. In contrast, our approach consciously uses an overly simplistic, but computationally light model, and relies heavily on tuning to find a feasible solution. 

To summarize, the main contributions of this paper are to (i) propose a computationally light-weight method to obtain high-quality solutions to chance-constrained problems, (ii) demonstrate the viability of this approach for OPF. The approach is applicable to both single and joint chance constrained problems and general uncertainty distributions.

The remainder of the paper is organized as follows: Section \ref{sec:formulation} presents the single and joint chance-constrained formulations of the optimal power flow problem. Section \ref{sec:reformulation} discusses the chance-constraint reformulation and the relationship to robust optimization, while Section \ref{sec:tuning} describes the tuning method. The case studies in Section \ref{sec:casestudy} demonstrate the efficiency and viability of the approach, while Section \ref{sec:conclusion} summarizes and discusses future directions. 

\section{Chance-Constrained Optimal Power Flow}
\label{sec:formulation}
We first present the formulation of the chance-constrained DC optimal power flow (CC-OPF) with both single and joint chance constraints, which is based on \cite{bienstock2014, vrakopoulou2013,  roaldSecurityConstrainedOptimal2015}.

\subsubsection*{Notation}
We represent a power system using an undirected graph $G = (\mathcal{N}, \mathcal{L})$, where $\mathcal{N}$ is the set of buses with $m = \abs{\mathcal{N}}$ and $\mathcal{L}$ is the set of lines with $l = \abs{\mathcal{L}}$. Without loss of generality, we assume that each bus has one generator, $g \in \mathcal{G} \subseteq{\real^m}$, one load $d \in \mathcal{D} \subseteq{\real^m}$, and one uncertainty source represented using random variable $\xi \in \Xi \subseteq \real^{m}$. Multiple (or zero) generators, loads, or uncertainty sources are handled by summation (or setting the respective elements to zero). We consider the DC linearized approximation of the power flow equations, where we assume voltages are constant at magnitude 1 per unit, angle differences are small, and the system is lossless.

\subsection{Modeling Considerations} 
\subsubsection{Uncertainty Modeling}
We consider uncertainty as a composite term, representing all fluctuations resulting from sources such as renewable generation and load uncertainty. The uncertainty fluctuations can be decomposed into $\xi = \mu + \omega$, where $\mu = \E[\xi] \in \real^m$ is the bias term and $\omega \in \real^m$ is the zero mean fluctuating component. The covariance of the fluctuations is denoted as $\Sigma_\xi = \text{Cov}[\xi] \in \real^{m \times m}$. Off-diagonals can be non-zero because we do not assume uncertainty sources are independent. We assume forecasts or estimates of $\mu$ and $\Sigma_u$ are available. For simplicity, we restrict ourselves to the case with unbiased forecasts, where $\mu = 0$.
 
\subsubsection{Power Balance and Generation Control}
Power systems operation requires power production and consumption to be balanced. The total power mismatch, $\Omega_\xi = \sum_{i = 1}^{m} \xi_i$, must be balanced by adjustments in controllable generation. We model this adjustment using an affine control policy based on actions of the automatic generation control (ACG), where $\Omega_\xi$ is divided amongst generators according to participation factors $\alpha \in [0,1]^{m}$ \cite{vrakopoulou2013}. We assume each generator contributes according to its maximum nominal output, i.e.,
\begin{align*}
 	\alpha_i = \frac{p_{G,i}^{\max}}{\sum_{j \in \mathcal{G}} p_{G,j}^{\max}} \qquad \forall i \in \mathcal{G}.
\end{align*}
The actual generation can be represented as
\begin{align*}
	\tilde{p}_{G,i}(\xi) = p_{G,i} - \alpha_i \Omega_\xi \qquad \forall i \in \mathcal{G},
\end{align*}
where $p_G$ denotes the scheduled generation. The total power balance is enforced with the following constraint
\begin{align*}
	\sum_{i \in \mathcal{N}} p_{G,i} - \alpha_i \Omega_\xi - d_i + \xi_i = 0.
\end{align*}
Because $\sum_{i \in \mathcal{G}}\alpha_i=1$ guarantees that any deviation $\xi$ is automatically balanced out with an equal adjustment in generation, it is sufficient to guarantee power balance for $\xi=0$, i.e.,
\begin{align*}
	\sum_{i \in \mathcal{N}} p_{G,i} - d_i = 0.
\end{align*}

\subsubsection{Power Flows}
We use a linear DC approximation to represent the power flow on lines $ij \in \mathcal{L}$ connecting buses $i, j \in \mathcal{N}$. We define $\mathbf{M} \in \real^{l \times m}$ to be the matrix of power transfer distributions factors (PTDFs) \cite{christie2000}, which relates the changes in active power flow to power injections at buses. Power flow on line $ij \in \mathcal{L}$ can be expressed as
\begin{align*}
	p_{ij} = \mathbf{M}_{(ij, \cdot)} (p_G - \alpha \Omega_\xi + \xi - d),
\end{align*}
where $\mathbf{M}_{(ij, \cdot)}$ is the row of $\mathbf{M}$ corresponding to line $ij$.

\subsubsection{Cost Function} The objective is to minimize the total generation cost of the scheduled generation $p_G$. Generation costs are modeled using a quadratic cost function, 
\begin{align}
	c(p_G) = \textstyle{\sum_{i \in \mathcal{G}}} \left(c_{2, i} p_{G,i}^2 + c_{1, i} p_{G, i} + c_{0,i}\right),
\end{align}
where the $c_{2,i}, c_{1,i}$ and $c_{i,0}$ are the quadratic, linear and constant cost coefficients corresponding to generator $i \in \mathcal{G}$. 

\subsection{Chance-Constrained OPF}
We consider both single and joint chance constraints, where $\epsilon$ denotes the acceptable violation probability. 

\subsubsection{Single Chance Constraints} 
\label{sec:scc-opf}
The OPF problem with single chance constraints (SCC-OPF) requires that each constraint is satisfied individually with separate acceptable violation probabilities. It can be formulated as
\begin{subequations}
\label{eqn:scc-opf}
\begin{align}
	&\min_{p_G} \enskip c(p_G) \\
	&\text{s.t.} \sum_{i \in \mathcal{N}} p_{G,i} - d_i = 0 \\
	& \mathbb{P}_{\xi} \!\left( p_{G,i} \!-\! \alpha_i \Omega_\xi \!\leq\! p_{G,i}^{\max} \right) \!\geq\! 1 \!-\! \epsilon, \quad \forall i \!\in\! \mathcal{G} \label{eqn:scc-opf-gen-max} \\
	& \mathbb{P}_{\xi} \!\left( p_{G,i} \!-\! \alpha_i \Omega_\xi \!\geq\! p_{G,i}^{\min} \right) \!\geq\! 1 \!-\! \epsilon, \quad \forall i \!\in\! \mathcal{G} \label{eqn:scc-opf-gen-min}\\
	& \mathbb{P}_{\xi} \!\left( \mathbf{M}_{(ij,\cdot)}(p_G \!-\! \alpha \Omega_\xi + \xi \!-\! d) \!\leq\! p_{ij}^{\max} \right) \!\geq\! 1 \!-\! \epsilon, 
	~\forall ij \!\in\! \mathcal{L} \\
	& \mathbb{P}_{\xi} \!\left( \mathbf{M}_{(ij,\cdot)}(p_G \!-\! \alpha \Omega_\xi + \xi \!-\! d) \!\geq\! - p_{ij}^{\max} \right) \!\geq\! 1 \!-\! \epsilon, 
	\forall ij \!\in\! \mathcal{L}
\end{align}
\end{subequations}
We define $\mathcal{C}$ as the set of all chance constraints, where $|\mathcal{C}| = 2m + 2l$. We note that the generator chance constraints \eqref{eqn:scc-opf-gen-max}, \eqref{eqn:scc-opf-gen-min} depend only on the total power mismatch $\Omega_\xi$, which is a scalar random variable. As a result, all generators will adjust their generation output up or down in \emph{perfect correlation}.

\subsubsection{Joint Chance Constraints} The formulation with a joint chance constraint requires all constraints to be simultaneously enforced with a single acceptable violation probability.
The joint chance-constrained OPF problem (JCC-OPF) can be formulated as
\begin{subequations}
\label{eqn:jcc-opf}
\begin{align}
	&\!\!\min_{p_G} \enskip c(p_G) \\
	&\!\text{s.t.}  \sum_{i \in \mathcal{N}} p_{G,i}  - d_i = 0 \\
	& \!\mathbb{P}_\xi \left( \begin{array}{l}
		\!\!p_{G,i} - \alpha_i \Omega_\xi \leq p_{G,i}^{\max}, ~~\forall i \in \mathcal{G} \\
		\!\!p_{G,i} - \alpha_i \Omega_\xi \geq p_{G,i}^{\min}, ~~\forall i \in \mathcal{G} \\
		\!\!\mathbf{M}_{(ij,\cdot)}(p_G - \alpha \Omega_\xi + \xi \!- \!d) \leq p_{ij}^{\max}, ~\forall ij \in \mathcal{L} \\
		\!\!\mathbf{M}_{(ij,\cdot)}(p_G - \alpha \Omega_\xi + \xi \!-\! d) \geq \!- p_{ij}^{\max}, \forall ij \in \mathcal{L}
 	\end{array} \!\!\right) 
 	\!\!\geq \!1 \!- \!\epsilon.
\end{align}
\end{subequations}

\section{Generalized chance constraint reformulation}
\label{sec:reformulation}
To become computationally tractable, the chance constraints in \eqref{eqn:scc-opf} and \eqref{eqn:jcc-opf} must be reformulated into detereministic constraints. This is challenging because (i) the probability term on the left hand side of the constraints is difficult to evaluate and (ii) the constraints often admit non-convex feasible sets. The goal of our approach is to use a simple reformulation and tune the parameters to achieve good performance. We therefore start from an analytical reformulation. 

\subsection{Individual Chance Constraint Reformulation}
In the case of individual chance constraints, we can obtain closed form deterministic reformulations under the assumption that the underlying distribution of $\xi$ is a Gaussian distribution. Specifically, assume that $\xi$ follows a multivariate Gaussian distribution parameterized by mean $\mu$ and covariance $\Sigma_\xi$. Consider the following chance constraint on generator $i \in \mathcal{G}$:
\begin{align}
    \label{eqn:scc-generator}
	\Pr \left( p_{G,i} - \alpha_i \Omega_\xi \leq p_{G,i}^{\max} \right) \geq 1 - \epsilon.
\end{align}
The deterministic reformulation is
\begin{align}
    p_{G,i} \leq p_{G,i}^{\max} - \Phi^{-1}(1 - \epsilon) ||\alpha_i \mathbbm{1}_{1,m} \Sigma_{\xi}^{1/2}||_2,
    \label{eqn:scc-reform-w-dist}
\end{align}
where $\Phi^{-1}(1 - \epsilon)$ is the inverse cumulative Gaussian evaluated at violation level $\epsilon$. This reformulation is tight, meaning that the chance constraint (\ref{eqn:scc-generator}) achieves a violation probability exactly equal to $\epsilon$ whenever reformulated constraint (\ref{eqn:scc-reform-w-dist}) is active (i.e., satisfied with equality). Moreover, each chance constraint only gives rise to \emph{one} reformulated constraint and this constraint is linear. By utilizing the analytical reformulation for all chance constraints in (\ref{eqn:scc-opf}), SCC-OPF becomes a linear program which is very efficiently solvable. 

However, this reformulation not only requires perfect knowledge of the distribution of $\xi$, but is also applicable only to the limited class of multivariate elliptical distributions, which includes distributions such as multivariate Gaussian, Student's $t$, and Cauchy distributions. Because these assumptions are very strong and may frequently not be satisfied in practical power system operations, this method may yield inaccurate (and possibly non-conservative) results. We can instead consider \textit{distributionally robust reformulations} where the distribution of $\xi$ is assumed to be uncertain within a family of distributions $\mathcal{P}$, known as an ambiguity set. Following \cite{roaldSecurityConstrainedOptimal2015}, we may generalize the $\Phi^{-1}(1 - \epsilon)$ term to $f^{-1}(1 - \epsilon)$, where $f^{-1}$ can be appropriately determined or bounded by invoking various inequalities depending on the assumed properties of $\mathcal{P}$.
These bounds are however usually only tight for one particular distribution $P \in \mathcal{P}$. As a result, the violation probabilities observed are generally less than $\epsilon$, leading to a solution that may be overly conservative and expensive.

\subsection{Robust Optimization Perspective}
We can consider the above analytical reformulation in a more general viewpoint by using connections to robust optimization. Specifically, the reformulated chance constraint (\ref{eqn:scc-opf}) can be interpreted as a robust constraint with an ellipsoidal uncertainty set $\mathcal{U}$ \cite{ben-talRobustOptimization2009}. To see this, consider the following robust constraint,
\begin{align}
    \label{eqn:scc-ro}
    \{ p_{G,i} - \alpha_i \Omega_\xi \leq p_{G,i}^{\max}: \quad \forall \xi \in \mathcal{U} \},
\end{align}
which requires the constraint to hold for all $\xi \in \mathcal{U}$. If we consider ellipsoidal uncertainty sets, $\mathcal{U}$ can be expressed as
\begin{align}
	\label{eqn:uncertainty-set}
    \mathcal{U} = \{ \xi: \norm{\xi}_2 \leq s \}.
\end{align}
Here $s$, known as the \textit{safety parameter} \cite{ben-talRobustOptimization2009,lam2019combating}, controls the size of the ellipsoid. The robust constraint (\ref{eqn:scc-ro}) with uncertainty set (\ref{eqn:uncertainty-set}) can be expressed as \cite{ben-talRobustOptimization2009}
\begin{align}
    \label{eqn:robust-const-reform}
	p_{G,i} \leq p_{G,i}^{\max} - s ||\alpha_i \mathbbm{1}_{1,m} \Sigma_{\xi}^{1/2}||_2.
\end{align}
Intuitively, by choosing $\mathcal{U}$ such that it contains sufficient probability mass, every solution satisfying this family of constraints (\ref{eqn:scc-ro}) will satisfy the original chance constraint with probability greater than $1 - \epsilon$.
In the case that $\xi \sim \mathcal{N}(\mu, \Sigma_\xi)$, by taking $s = \Phi^{-1}(1 - \epsilon)$, constraint (\ref{eqn:robust-const-reform}) is exactly equivalent to (\ref{eqn:scc-reform-w-dist}). A larger $s$ value will correspond to a higher probability that the inner constraint of the chance constraint will be satisfied because the uncertainty set contains more probability mass.

\subsection{SCC-OPF Reformulation} By replacing the $\Phi^{-1}(1 - \epsilon)$ term in (\ref{eqn:scc-reform-w-dist}) with the safety parameter $s$, we generalize the analytical reformulation such that it does not rely on any distributional assumptions on $\xi$. Applying this formulation to all chance constraints in SCC-OPF, we obtain
\begin{subequations}
\label{eqn:analytical-reform-s}
\begin{align}
	\min_{p_G} \enskip & c(p_G) \\
	\text{s.t.} & \sum_{i \in \mathcal{N}} p_{G,i} - d_i + = 0 \\
	& p_{G,i} \leq p_{G,i}^{\max} - s ||\alpha_i \mathbbm{1}_{1,m} \Sigma_{\xi}^{1/2}||_2, \forall i \in \mathcal{G} \\
	& p_{G,i} \geq p_{G,i}^{\min} + s ||\alpha_i \mathbbm{1}_{1,m} \Sigma_{\xi}^{1/2}||_2, \forall i \in \mathcal{G} \\
	& \mathbf{M}_{(ij, \cdot)}(p_G - d) \leq p_{ij}^{\max} \nonumber \\
	& \qquad - s || \mathbf{M}_{(ij,\cdot)} (I - \alpha \mathbbm{1}_{1,m}) \Sigma_\xi^{1/2} ||_2, \forall ij \in \mathcal{L} \\
	& \mathbf{M}_{(ij, \cdot)}(p_G - d) \geq - p_{ij}^{\max} \nonumber \\ 
	& \qquad + s || \mathbf{M}_{(ij,\cdot)} (I - \alpha \mathbbm{1}_{1,m}) \Sigma_\xi^{1/2} ||_2, \forall ij \in \mathcal{L}.
\end{align}
\end{subequations}
We observe that if we can determine an appropriate $s$ value, we will obtain a solution that is tight for at least one constraint in the SCC-OPF with a pre-specified violation probability. This forms the basis of our proposed tuning method.

\subsection{JCC-OPF Reformulation}
For joint chance constraints, there is no analogous analytical reformulation. It is possible to obtain an upper bound on the chance constraint using Boole's inequality 
\cite{boole1854investigation} 
by separating the joint constraint into $k$ individual constraints, each with violation level $\epsilon/k$. Unfortunately, this method typically results overly conservative solutions, particularly when only a small number of constraints experience violations or when some constraints are perfectly correlated.
Alternative approaches that achieve tighter bounds have been proposed in, e.g., \cite{bakerEfficientRelaxationsJoint2016}, but the methods remain conservative.

We instead propose that by accurately tuning the safety parameter $s$, the generalized reformulation for single chance constraints (\ref{eqn:analytical-reform-s}) can be carried over to joint chance constraints.

\section{Chance Constraint Tuning}
\label{sec:tuning}
Our goal is to tune the safety parameter $s$ such that the solution to the corresponding analytical reformulation (\ref{eqn:analytical-reform-s}) exactly satisfies our desired joint or single violation probability, $\epsilon_{\text{des}}$.
Due to the monotonic relationship between $s$ and the level of violation (i.e., a larger $s$ leads to tighter constraints, a more conservative solution, and lower violation probabilities), our proposed method uses a bisection search to determine the value of $s$.
A similar approach is used in \cite{pena2019dc}.

\subsection{Bisection Search}
\begin{itemize}
	\item[0)] 
	\noindent\emph{Initialization:} We first set the iteration count to $k=0$ and determine suitable upper and lower bounds for $s$:\\ 
	For the \emph{lower bound}, we use $s_{\min} = 0$, which corresponds to the case where no uncertainty is considered in the generation and line flow constraints. Because $s$ has an inverse monotonic correlation with the observed empirical violation probability, $\epsilon_{\text{obs}}$, we expect that the solution to (\ref{eqn:analytical-reform-s}) corresponding to $s = 0$ will have the highest possible $\epsilon_{\text{obs}}$ for any $s \geq 0$. \\
	{
	For the \emph{upper bound}, we use the Cantelli inequality \cite{roaldSecurityConstrainedOptimal2015} since we have estimates of $\mu$ and $\Sigma_{\xi}$.\\[+2pt]
	For SCC-OPF, we set $s_{\max}~=~\sqrt{(1 - \epsilon_{\text{des}})/\epsilon_{\text{des}}}$.\\[+2pt]
	For JCC-OPF, we use this in conjunction with Boole's inequality, and set $s_{\max} = \sqrt{(1 - \frac{\epsilon_{\text{desired}}}{|\mathcal{C}|})/\frac{\epsilon_{\text{des}}}{|\mathcal{C}|}}$. \\[+2pt]
	These choices of $s_{\max}$ guarantee that $\epsilon_{\text{obs}}\leq \epsilon_{\text{des}}$.\\[-10pt]}
		
	\item[1)]
	\emph{Solve OPF:} We increase the iteration count to $k=k+1$ and define $s_k = (s_{\max} - s_{\min})/2 + s_{\min}$.
	We then solve (\ref{eqn:analytical-reform-s}) with $s=s_k$ to obtain a candidate solution. If (\ref{eqn:analytical-reform-s}) is infeasible, we set $s = s_{\max}$ and repeat this step.\\[-10pt]
	\item[2)]
	\emph{Evaluate violations:} Using samples $\{\xi^{(1)}, \dots, \xi^{(N)} \}$, we evaluate either the observed worst-case single empirical violation probability, $\epsilon_{\text{obs, joint}}$, or observed joint violation probability, $\epsilon_{\text{obs, single}}$, depending on whether we are solving JCC-OPF or SCC-OPF. The exact definitions of these values are described in Section \ref{sec:empirical}. 
	\item[3)] 
	\emph{Update $s$:} We update $s$ based on observed violations. \\ If $\epsilon_{\text{obs}} < \epsilon_{\text{des}}$, the current solution is too conservative, meaning $s_k$ is larger than the value we are looking for. To decrease $s$ in the next iteration, we set $s_{\max}=s_k$. \\ If $\epsilon_{\text{obs}} > \epsilon_{\text{des}}$, the current $s_k$ is too small. To increase $s$ in the next iteration, we set $s_{\min}=s$.\\[-10pt]
	\item[4)] 
	\emph{Check convergence:} This process is repeated, using the same set of evaluation samples $\{\xi^{(1)}, \dots, \xi^{(N)} \}$, 
	until $|\epsilon_{\text{obs}} - \epsilon_{\text{des}}| \leq \gamma$, where $\gamma$ is our tolerance value. Convergence of the bisection search is guaranteed to occur in $\lfloor \log_2\big( (s_{\max} - s_{\min})/\gamma \big) \rfloor$ iterations.
\end{itemize}

\subsection{Evaluation of Empirical Violation Probability}
\label{sec:empirical}
\subsubsection{Single Chance Constraint}
We first check whether the inner inequality of each chance constraint $\ell \in \mathcal{C}$ holds for each sample $\{\xi^{(1)}, \dots, \xi^{(N)} \}$. In particular, taking constraint (\ref{eqn:scc-opf-gen-max}) as an example, for sample $\xi^{(n)}$, we count whether a violation has occurred using
\begin{align}
    V_{\ell}(\xi^{(n)}) = 
    \begin{cases}
    	1 & \text{ if } p_{G,i} - \alpha_i \Omega_{\xi^{(n)}} > p_{G,i}^{\max} \\
    	0 & \text{ otherwise. }
    \end{cases}
\end{align}
The observed empirical violation probability is
\begin{align}
    \textstyle
	\epsilon_{\text{obs}}^{\ell} = \frac{1}{N} \sum_{n = 1}^{N} V_{\ell} ({\xi^{(n)}}).
\end{align}
Calculating this for all $\ell \in \mathcal{C}$, the worst case empirical violation probability is found by taking the maximum
\begin{align}
	\epsilon_{\text{obs, single}} = \max_{\ell \in \mathcal{C}} \epsilon_{\text{obs}}^{\ell}.
\end{align}

\subsubsection{Joint Chance Constraint}
The empirical violation probability is determined similarly, but we consider a violation to have occurred if any one constraint is violated. The observed empirical violation probability is
\begin{align}
    \textstyle
	\epsilon_{\text{obs, joint}} = \frac{1}{N} \sum_{n = 1}^{N} \max_{\ell \in \mathcal{C}} V_{\ell}.
\end{align}

\begin{table*}
\centering
	\begin{tabular}{|| c | c | c c c c c c c c c ||} 
	\hline
	\multicolumn{2}{|| c}{} & $\epsilon_{\text{des}}$ & iterations	& costs &	$s$ & $s_\text{true}$ & $\epsilon_{\text{obs, single}}$ & $\epsilon_{\text{oos, single}}$ & $\epsilon_{\text{obs, joint}}$ & $\epsilon_{\text{oos, joint}}$ \\ [0.5ex] 
	\hline \hline
	\multirow{6}{4em}{Single constraints} &
	\multirow{3}{4em}{Gaussian}
	& 0.10		& 10.7  & 42201.6 & 1.3012 & 1.2816 & 0.1000	& 0.0976	& 0.2960 	& 0.2947 \\
	& & 0.05	& 9.4 	& 42376.1 & 1.6676 & 1.6449 & 0.0501	& 0.0483 	& 0.1585	& 0.1577 \\
	& & 0.01	& 9.6 	& 42709.0 & 2.3624 & 2.3263 & 0.0100 & 0.0095 	& 0.0338	& 0.0338 \\
	\cline{2-11} 
	& \multirow{3}{4em}{Non-Gaussian}
	& 0.10		& 10.4  & 42799.6 & 1.3376	& - & 0.1001	& 0.1007	& 0.3031	& 0.3044  \\
	& & 0.05	& 9.6 	& 43105.4 & 1.6677 & - & 0.0501 	& 0.0495	& 0.1597	& 0.1609 \\
	& & 0.01	& 8.9  	& 43680.4 & 2.2844 & -	& 0.0100 	& 0.0095	& 0.0274 	& 0.0275 \\
	\hline \hline
	\multirow{6}{4em}{Joint constraints} &
	\multirow{3}{4em}{Gaussian}
	& 0.10	    & 15.6 	& 42485.6 & 1.8971 	& -	& 0.0307	& 0.0296 	& 0.1001	& 0.1001 \\
	& & 0.05 	& 14.5 	& 42632.9 & 2.2054	& - & 0.0149	& 0.0141 	& 0.0501 	& 0.0500 \\
	& & 0.01	& 12.1 	& 42918.5 & 2.8014 	& -	& 0.0032	& 0.0027 	& 0.0100 	& 0.0100 \\
	\cline{2-11} 
	& \multirow{3}{4em}{Non-Gaussian}
	& 0.10	    & 14.4 	& 43284.4	& 1.8585	& - & 0.0316	& 0.0307 	& 0.1001	& 0.1000 \\
	& & 0.05 	& 14.6 	& 43507.0	& 2.1008    & - & 0.0167 	& 0.0161 	& 0.0500	& 0.0501 \\
	& & 0.01	& 13.3 	& 43924.0	& 2.5538    & - & 0.0042 	& 0.0038 	& 0.0100	& 0.0101 \\
	\hline
	\end{tabular}
	\caption{\small Results for the bisection tuning method for joint and single chance-constrained OPF (average values for 20 replications).}
\label{table:jcc-scc-gaussian-non-gaussian}
\end{table*}

\section{Case Studies}
\label{sec:casestudy}

\subsection{Test system}
We evaluate our method on the IEEE RTS96 24-bus system \cite{griggIEEEReliabilityTest1999}, with the following modifications: (i) line capacities are reduced to 70\%, (ii) minimum output is set to 0 on all generators, and (iii) maximum output is doubled on all generators. On buses 8 and 15, we add uncertainty sources to represent variations in load and renewable energy. For all experiments, we run 20 replications with tolerance $\gamma =10^{-4}$. We use 10,000 samples of $\xi$ during the tuning process and 100,000 samples for the out-of-sample evaluation.

\begin{figure*}
    \centering
    \includegraphics[width=.95\linewidth]{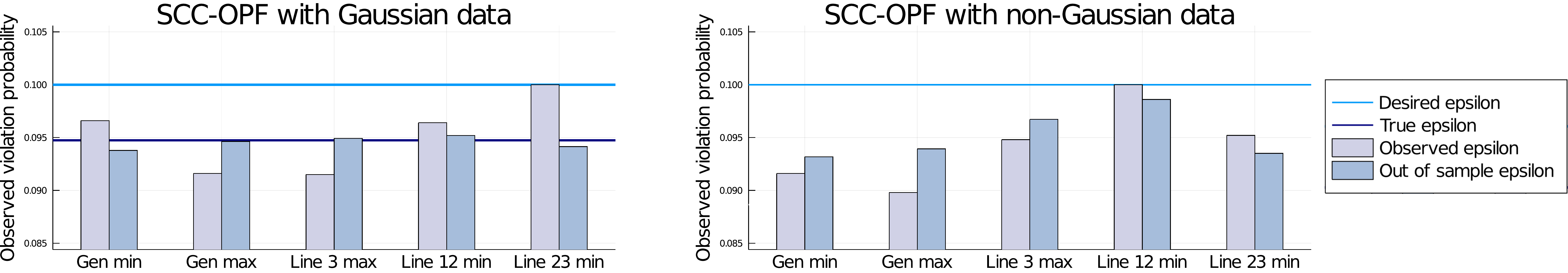}
    \caption{\small Empirical violation probability for the active constraints in the SCC-OPF problem, with Gaussian (left) and non-Gaussian (right) data. The active constraints includes several generator maximum and minimum limits (represented by one bar because they are perfectly correlated and hence have the same violation probability), as well as three line limits. The blue and grey bars show the empirical violation probability of each constraint, as observed in the tuning and out-of-sample evaluation respectively. The light blue line shows $\epsilon_{\text{des}}$, while the dark blue line is $\epsilon_{s} = 1 - \Phi(s)$.}%
    \label{fig:scc-opf}
\end{figure*}

\begin{figure*}
    \centering
    \includegraphics[width=.95\linewidth]{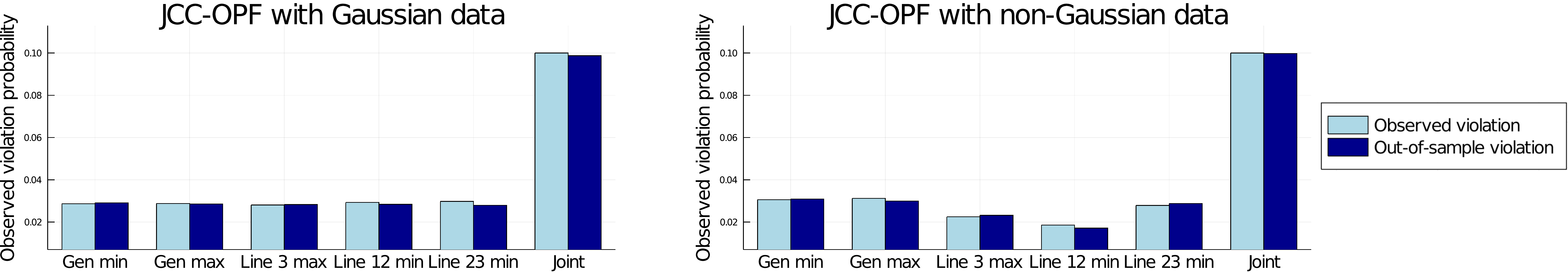}
    
    \caption{\small Empirical violation probability for individual active constraints in the JCC-OPF problem, as well as the joint violation probability. We show results for both Gaussian data (left) and non-Gaussian data (right). The grey and blue bars represent the empirical violation probability of each constraints, as observed in the tuning and out-of-sample evaluation, respectively.}
    \label{fig:joint}
\end{figure*}

\subsection{Proof-of-concept: SCC-OPF with Gaussian data}
\label{section:scc-gaussian}
To demonstrate the accuracy and efficiency of our method, we perform a proof of concept study by evaluating performance on the SCC-OPF problem using Gaussian uncertainty data. In this case, we can directly calculate the true safety parameter value by evaluating $s_{\text{true}} = \Phi^{-1}(1 - \epsilon_{\text{des}})$. We test whether our method identifies a parameter $s\approx s_{\text{true}}$ by tuning the observed violation probability $\epsilon_{\text{obs,single}}$ to match $\epsilon_{\text{des}}$.

Data is generated from a multivariate Gaussian distribution with zero mean, standard deviations of 9.4 MW (bus 8) and 13.1 MW (bus 15), and a correlation coefficient of $\rho = 0.2$. We consider desired epsilon values $\epsilon_{\text{des}} = \{0.1, 0.05, 0.01\}$, which correspond to true safety parameter values $ s_{\text{true}} = \{1.2816, 1.6449, 2.3263 \}$. 

The top section of Table \ref{table:jcc-scc-gaussian-non-gaussian} shows the the cost of generation, resulting safety parameter values, resulting in-sample single and joint empirical violation probabilities, and out-of-sample joint and single violation probabilities. All number represent average values across the 20 replications. 
We observe that the number of iterations is about $10$ in all cases, and that the method terminates with $\epsilon_{\text{obs,single}}$ very close to $\epsilon_{\text{des}}$. Furthermore, using $10,000$ samples to evaluate $\epsilon_{\text{obs,single}}$ in the tuning step is sufficient to achieve out-of-sample violation probabilities $\epsilon_{\text{oos,single}}$ that are very close to the desired values. Finally, we conclude that the average resulting $s$ is close to $s_{\text{true}}$, verifying the correctness and accuracy of our method.

An interesting, but subtle observation, is that the $s$ values obtained with the algorithm are typically always slightly \emph{greater} (i.e., more conservative) than $s_{\text{true}}$ for every $\epsilon_{\text{desired}}$.
This is true not only for the average $s$ in Table \ref{table:jcc-scc-gaussian-non-gaussian}, but for each of the 20 runs. To explain this behavior we look at the results of a single run in more detail, as shown in Figure \ref{fig:scc-opf} (right). 
We see that the tuning algorithm chooses the value of $s$ to ensure that the worst case single violation probability $\epsilon_{\text{obs,single}} \leq 0.1$, leading to a conservative result. In fact, because we are using Gaussian data, we can calculate the true violation probability corresponding to the tuned value of $s$, $\epsilon_{s} = 1 - \Phi(s)\approx0.095$. We observe that the out-of-sample violation probabilities $\epsilon_{\text{o.o.s.,single}}$ (grey bars) are all close to this true violations probability. 

\subsection{SCC-OPF with non-Gaussian data}
We next study how the bisection method performs for SCC-OPF using non-Gaussian data. To generate a non-Gaussian data set, we mix the following distributions: (i) a zero-mean multivariate Gaussian with standard deviations of $7$ MW and $14$ MW and a correlation coefficient of $0.5$, (ii) a zero-mean multivariate Gaussian with standard deviations of both $6$ MW and a correlation coefficient of $0.1$, and (iii) a uniform distribution on the interval $[-30,30]$.  We consider desired violation probabilities $\epsilon_{\text{desired}} = \{0.1, 0.05, 0.01 \}$. The results are reported in the second section of Table \ref{table:jcc-scc-gaussian-non-gaussian}.

We observe that the number of iterations is not impacted by the type of uncertainty data, and that the algorithm also terminates with $\epsilon_{\text{obs,single}}\approx\epsilon_{\text{desired}}$ for the non-Gaussian data. Again, 10,000 samples in the tuning algorithm is sufficient to obtain an $\epsilon_{\text{o.o.s.,single}}$ that is very close to the desired values. 

Figure \ref{fig:scc-opf} shows the result for a single run in more detail. We observe larger variations in $\epsilon_{\text{obs, single}}$ and $\epsilon_{\text{oos, single}}$ amongst the constraints compared to using Gaussian data. Because we do not have the Gaussian distribution assumption, it is no longer true that the same safety parameter, $s$, will lead to the same true violation probability for all constraints. 

\subsection{JCC-OPF with Gaussian and non-Gaussian data}
We finally investigate the behavior of the algorithm for JCC-OPF. We use the same case study set up and the same Gaussian and non-Gaussian distributions as above, but tune the value of $s$ to achieve a desired \emph{joint} violation probability. The results are shown in the lower half of Table \ref{table:jcc-scc-gaussian-non-gaussian}. 

We observe that the number of iterations remains similar, but is slightly higher compared to the number of iterations for the single chance constraints. As for the single chance constraints, the algorithm manages to determine an $s$ value which meets the joint violation probability $\epsilon_{\text{obs,joint}}$ exactly. However, the out-of-sample results show, on average, violation probabilities that are slightly too high for most cases, indicating that the results are no longer conservative. 

To gain some more insight into the solutions for the JCC-OPF, in Figure \ref{fig:joint} we plot the individual and joint constraint violation levels for JCC-OPF with Gaussian and non-Gaussian data. We again observe that the violation probability is spread more evenly in the case of Gaussian data, while the non-Gaussian data lead to larger variations despite all constraints sharing the same $s$. For both data types, the sum of the violation probabilities of the individual chance constraints exceeds the joint violation probability, i.e., $\sum_{\ell \in \mathcal{C}} \epsilon_{\text{obs}}^{\ell} \geq \epsilon_{\text{obs,joint}}$. This indicates that our algorithm is able to account for correlation between different constraints. 

\section{Conclusion and Future Work}
\label{sec:conclusion}
In summary, we propose a method that uses samples of uncertainty realizations for a posteriori tuning of chance constraints. The benefits of our approach include low computational overhead, applicability to both single and joint chance constraints, and ability to handle general uncertainty distributions.

For future work, we would like to obtain theoretical guarantees on sample complexity and error bounds. Moreover, we would like to separately tune the individual chance constraints, which would give us additional degrees of freedom and therefore possibly more optimal solutions. This extension is non-trivial due to the large number of constraints and resulting high-dimensional tuning challenge. Finally, we note that it may be possible to extend this method to solve AC-OPF. However, the non-linearity poses challenges, as an inverse monotonic relationship between the safety parameter $s$ and the observed empirical violation probability no longer necessarily holds.

\bibliographystyle{IEEEtran}
\bibliography{ref}

\begin{thebibliography}{10}
\providecommand{\url}[1]{#1}
\csname url@samestyle\endcsname
\providecommand{\newblock}{\relax}
\providecommand{\bibinfo}[2]{#2}
\providecommand{\BIBentrySTDinterwordspacing}{\spaceskip=0pt\relax}
\providecommand{\BIBentryALTinterwordstretchfactor}{4}
\providecommand{\BIBentryALTinterwordspacing}{\spaceskip=\fontdimen2\font plus
\BIBentryALTinterwordstretchfactor\fontdimen3\font minus
  \fontdimen4\font\relax}
\providecommand{\BIBforeignlanguage}[2]{{%
\expandafter\ifx\csname l@#1\endcsname\relax
\typeout{** WARNING: IEEEtran.bst: No hyphenation pattern has been}%
\typeout{** loaded for the language `#1'. Using the pattern for}%
\typeout{** the default language instead.}%
\else
\language=\csname l@#1\endcsname
\fi
#2}}
\providecommand{\BIBdecl}{\relax}
\BIBdecl

\bibitem{abbaspourtorbati2016}
F.~Abbaspourtorbati and M.~Zima, ``{The Swiss Reserve Market: Stochastic
  Programming in Practice},'' \emph{{IEEE Trans. Power Systems}}, vol.~31,
  no.~2, pp. 1188--1194, March 2016.

\bibitem{roaldAnalyticalReformulationSecurity2013}
L.~Roald, F.~Oldewurtel, T.~Krause, and G.~Andersson,
  ``\BIBforeignlanguage{en}{Analytical reformulation of security constrained
  optimal power flow with probabilistic constraints},'' in
  \emph{\BIBforeignlanguage{en}{IEEE PowerTech Conference}}, Jun. 2013.

\bibitem{bienstock2014}
D.~Bienstock, M.~Chertkov, and S.~Harnett, ``Chance-constrained optimal power
  flow: Risk-aware network control under uncertainty,'' \emph{SIAM Review},
  vol.~56, no.~3, pp. 461--495, 2014.

\bibitem{calafioreScenarioApproachRobust2006}
G.~Calafiore and M.~Campi, ``\BIBforeignlanguage{en}{The {{Scenario Approach}}
  to {{Robust Control Design}}},'' \emph{\BIBforeignlanguage{en}{IEEE Trans.
  Automatic Control}}, vol.~51, no.~5, pp. 742--753, May 2006.

\bibitem{vrakopoulou2013}
M.~Vrakopoulou, K.~Margellos, J.~Lygeros, and G.~Andersson, ``{A Probabilistic
  Framework for Reserve Scheduling and N-1 Security Assessment of Systems With
  High Wind Power Penetration},'' \emph{IEEE Trans. Power Systems}, vol.~28,
  no.~4, pp. 3885--3896, 2013.

\bibitem{summers2014}
T.~Summers, J.~Warrington, M.~Morari, and J.~Lygeros, ``Stochastic optimal
  power flow based on convex approximations of chance constraints,'' in
  \emph{Power System Computation Conference, 2014}, Aug 2014, pp. 1--7.

\bibitem{roaldSecurityConstrainedOptimal2015}
L.~Roald, F.~Oldewurtel, B.~Van~Parys, and G.~Andersson,
  ``\BIBforeignlanguage{en}{Security {{Constrained Optimal Power Flow}} with
  {{Distributionally Robust Chance Constraints}}},''
  \emph{\BIBforeignlanguage{en}{arXiv:1508.06061 [math]}}, Aug. 2015.

\bibitem{lubin2015}
M.~Lubin, Y.~Dvorkin, and S.~Backhaus, ``A robust approach to chance
  constrained optimal power flow with renewable generation,'' \emph{IEEE Trans.
  Power Systems}, vol.~31, no.~5, pp. 3840--3849, 2015.

\bibitem{li2019distributionally}
B.~Li, R.~Jiang, and J.~L. Mathieu, ``Distributionally robust
  chance-constrained optimal power flow assuming unimodal distributions with
  misspecified modes,'' \emph{IEEE Trans. Control of Network Systems}, vol.~6,
  no.~3, pp. 1223--1234, 2019.

\bibitem{luedtkeSampleApproximationApproach}
J.~Luedtke and S.~Ahmed, ``\BIBforeignlanguage{en}{A {{Sample Approximation
  Approach}} for {{Optimization}} with {{Probabilistic Constraints}}
  {${_\ast}$}},'' p.~23.

\bibitem{pena2019dc}
A.~Pena-Ordieres, D.~Molzahn, L.~Roald, and A.~Waechter, ``Dc optimal power
  flow with joint chance constraints,'' \emph{arXiv preprint arXiv:1911.12439},
  2019.

\bibitem{lam2019combating}
H.~Lam and H.~Qian, ``Combating conservativeness in data-driven optimization
  under uncertainty: A solution path approach,'' \emph{arXiv preprint
  arXiv:1909.06477}, 2019.

\bibitem{mezghani2019stochastic}
I.~Mezghani, S.~Misra, and D.~Deka, ``Stochastic ac optimal power flow: A
  data-driven approach,'' \emph{arXiv preprint arXiv:1910.09144}, 2019.

\bibitem{oldewurtelAdaptivelyConstrainedStochastic2013}
F.~Oldewurtel, D.~Sturzenegger, P.~M. Esfahani, G.~Andersson, M.~Morari, and
  J.~Lygeros, ``Adaptively constrained stochastic model predictive control for
  closed-loop constraint satisfaction,'' in \emph{IEEE American Control
  Conference}, 2013, pp. 4674--4681.

\bibitem{oldewurtelAdaptivelyConstrainedStochastic2015}
F.~Oldewurtel, L.~Roald, G.~Andersson, and C.~Tomlin, ``Adaptively constrained
  stochastic model predictive control applied to security constrained optimal
  power flow,'' in \emph{IEEE American Control Conference}, 2015, pp. 931--936.

\bibitem{christie2000}
R.~Christie, B.~F. Wollenberg, and I.~Wangensteen, ``Transmission management in
  the deregulated environment,'' \emph{Proceedings of the IEEE}, vol.~88,
  no.~2, pp. 170--195, Feb. 2000.

\bibitem{ben-talRobustOptimization2009}
A.~Ben-Tal, L.~El~Ghaoui, and A.~Nemirovski, \emph{Robust optimization}.\hskip
  1em plus 0.5em minus 0.4em\relax Princeton University Press, 2009, vol.~28.

\bibitem{boole1854investigation}
G.~Boole, \emph{An investigation of the laws of thought: on which are founded
  the mathematical theories of logic and probabilities}.\hskip 1em plus 0.5em
  minus 0.4em\relax Dover Publications, 1854.

\bibitem{bakerEfficientRelaxationsJoint2016}
K.~Baker and B.~Toomey, ``Efficient relaxations for joint chance constrained ac
  optimal power flow,'' \emph{Electric Power Systems Research}, vol. 148, pp.
  230 -- 236, 2017.

\bibitem{griggIEEEReliabilityTest1999}
C.~Grigg \emph{et~al.}, ``\BIBforeignlanguage{en}{The {{IEEE Reliability Test
  System}}-1996},'' \emph{\BIBforeignlanguage{en}{IEEE Trans. Power Systems}},
  vol.~14, no.~3, pp. 1010--1020, Aug./1999.

\end{thebibliography}

\end{document}